\documentclass[12pt,draft,twoside]{article}
\usepackage[cp1250]{inputenc}
\usepackage{latexsym,amssymb,amsmath}
\title{On rectifiable spaces and its algebraical equivalents, topological
algebraic systems and Mal'cev algebras (continuation)}
\date{}
\author{N. I. Sandu}
\begin{document}
\maketitle

This paper is a natural continuation of paper "On rectifiable
spaces and its algebraical equivalents, topological algebraic
systems and Mal'cev algebras"  published in:
arxiv.org/abs/1309.4572. Thus we justify the need to pre\-sent the
entire material in an unified manner. This paper is the
continuation of Section 6 from the first paper.  It specifies and
corrects the roughest mistakes,  incorrect statements and nonsense
of the introduced concepts, which are available in numerous papers
on topological algebraic systems, basically in papers of
Academician Choban M. M. and his disciples. \vspace*{0.1cm}

Remark 4 for a topological algebraic system $\mathcal{A}$ with
defining space $X$  defines the notion of free topology. Now,
similarly to  \cite[pag. 182]{Mal3} we will show how this
topology, irrespective of the topology given in $\mathcal{A}$ ,
can be obtained from the topology of space $X$.

A set $U$ of elements of $\mathcal{A}$ will be called $X_0$-open
if for each polynomial $f$, in particular $f(x) = x$, and for each
system $x_1, \ldots, x_n$ of elements from $X$, for which

$f(x_1, \ldots, x_n) \in U$ there exist a neighborhood these
elements in $X$ that $f(U_1, \ldots, U_n) \subseteq U$. Obviously,
the intersection of finite system and the reunion of arbitrary
system of $X_0$-open sets is a $X_0$-open sets. Hence the
$X_0$-open sets define a topology on $\mathcal{A}$, which we will
call $X_0$-topology.

 For transfinite index $\lambda$ we define
$X_{\lambda}$-topology by induction. In particular, let for some
transfinite  $\lambda$ the topology $X_{\lambda}$ be defined. A
set $U$ from $\mathcal{A}$ will be called $X_{\lambda + 1}$-open
if for all basic operation $f(x_1, \ldots, x_n)$ and all $a_1,
\ldots, a_n$ from $\mathcal{A}$, satisfying the condition $f(a_1,
\ldots, a_n) \in U$, there are such $X_{\lambda}$-neighborhood
$U_1, \ldots, U_n$ of points $a_1, \ldots, a_n$ in $\mathcal{A}$
that $f(U_1, \ldots, U_n) \subseteq U$. Further, if the topologies
$X_{\lambda}$ are defined for all $\lambda$ less than the limit
$\alpha$, then the set $U$ will called $X_{\alpha}$-open, if it is
$X_{\lambda}$-open for all $\lambda < \alpha$.

The topologies $X_{\lambda}$ decrease monotonically with the
increase of $\lambda$ and all of them contain every admissible
topology of $\mathcal{A}$, i.e. not changing the topologies $X$
and making continuous the basic operations.

Literally repeating the corresponding judgments from \cite[pag.
182]{Mal3} it is proved.\vspace*{0.1cm}

8a). \textit{Let $\tau$ be the least transfinite, for which
$X_{\tau} = X_{\tau + 1}$. Then $X_{\tau}$ is the sought topology
of algebra $\mathcal{A}$, free  regarding to $X$.}\vspace*{0.1cm}

According to item 1e) for topological algebraic systems
($\Omega$-algebras) we ignore the property of continuous signature
and according to Proposition 1  we will consider that the basic
operations of topological algebraic systems ($\Omega$-algebras)
have a finite arity. For the analysis we quote the item $1e_1)$
 (see \cite{Chob3}) in a modified version, i.e.
we will consider that the basic operations of topological
algebraic systems ($\Omega$-algebras) have a finite
arity.\vspace*{0.1cm}

9a). \textmd{Fix a class of topological systems $\frak{K}$.}
\textmd{Let's allocate the following properties, which can have
the class $\frak{K}$.}

\textmd{Condition 1m. Closed with respect to subsystems.}

\textmd{Condition 2m. Closed with respect to Tychonoff products.}

\textmd{Condition 3m. All objects from $\frak{K}$ satisfy a given
set of  topological, algebraical properties $Q$ (the set $Q$ can
be $\{\emptyset\}$ or the separation axioms $T_0, T_1, T_2, \break
T_3$, or the requirement of a regularity, or the requirement of a
completely regularity and others).}

\textmd{Condition 4m. If $(G, \tau) \in \frak{K}$, $\tau_d$ is the
discrete topology and  $(G, \tau_d)$ is a topological system of
continuous signature, then $(G, \tau_d) \in \frak{K}$.}

\textmd{Condition 5m. Closed with respect to continuous
homomorphic images with property $Q$.}

\textmd{Condition 6m. If $(G, \tau_1) \in \frak K$ and with
respect to topology $\tau$ the pair $(G, \tau)$ is a topological
system with properties $Q$ of continuous signature $S$, then $(G,
\tau) \in \frak K$.}

9b). \textmd{Definition 4. A class $\mathcal{K}$ is called:}

\textmd{1. Topological $Q$-quasivariety, if  the conditions  1m -
4m are satisfied.}

\textmd{2. Topological $Q$-prequasivariety, if  the conditions 1m
- 5m are satisfied.}

\textmd{3. Topological $Q$-variety, if  the conditions  1m - 6m
are satisfied.}

\textmd{4. Topological complete $Q$-quasivariety, if  the
conditions 1m - 4m and 6m are satisfied.}

9c). \textmd{Definition 5. A class $\frak{K}$ of algebraical
systems of signature $S = \Omega \cup P$ is called quasivariety,
if}

\textmd{1. The class $\frak{K}$ is closed with respect to
subsystems;}

\textmd{2. The class $\frak{K}$ is closed with respect to
Cartesian product.} $1e_1)$.

9d). Recall, that an algebraic system $\mathcal{A}_e$ of signature
$\Omega$ is called unitary if it is from one element $e \in
\mathcal{A}$, all its basic predicates have the value true and
$F(e, \ldots, e) = e$ for any basic operation $F \in \Omega$. For
any signature $\Omega$ there exists an unique to an accuracy of
homeomorphism unitary system $\mathcal{A}_e$.

9e). We consider the condition: 7m. $\frak{K}$ contains an unitary
system $\mathcal{A}_e$.

9f). Obviously, condition 7m from item 9e) does not contradict
conditions 1m -- 6m and does not influence the topologies of
topological algebraic systems of class $\frak K$ from item 9a).

9g). The definition 5 from item 9c) is incorrect, does not
correspond to classical definition of quasivariety from
Proposition 3. To class $\frak{K}$ of topological algebraic
systems with conditions 1m, 2m from item 9c) we add the condition
7m from item 9e). Then, by Lemma 1  the class $\frak{K}$ becomes
an almost quasivariety and by Theorem 7 $\frak{K}$ contains free
topological algebraic systems $\mathcal{F}_m$ of any given rank $m
\geq 1$.

9h). Suppose that the almost quasivariety $\frak{K}$ from item 9g)
meets the condition 4m from item 9a) and let $\mathcal{F}_m$ be a
free topological system of $\frak{K}$. By condition 4m $\frak{K}$
contain free discrete algebraic systems $(\mathcal{F}_m, \tau_d)$.
But by item $10m_6)$ such a topological free $E$-algebra $F(X, K)$
of quasivariety of topological $E$-algebras with discrete space
$X$ do not exist.

$9h_1)$. In particular, a topological non-discrete free
$E$-algebra $F(X, K)$ does not exist if $K$ is an almost
quasivariety and $X$ is a finite space since any topological
finite space is discrete.

9i). From items 9h) and 9f) it follows that all definitions from
item 9b) and item 9c) according to item 9g) are totally senseless.
We mentioned that besides paper \cite{Chob3} this non-sense is
investigated almost in almost all analyzed papers: \cite{Chob7},
\cite{Chir1}, \cite{Chir2}, \cite{Chob2}, \cite{Chob9},
\cite{Chob19}, \cite{Chob32}, \cite{Chob33}, \cite{Chob191},
\cite{Chob192}, \cite{Chob12}, \cite{Chob6}, \cite{ChobChir44}
\cite{Dum}, \cite{DC},   \cite{Chob28}. \vspace*{0.1cm}

$9i_1)$.All results from \cite{Chob20} are false, \cite{CD2},
which consider classes of topological algebras with conditions 1m
-- 4m from item 9a).\vspace*{0.1cm}

10a). The paper \cite{ChobChob} is published in ROMAI Journal and
can be accessed freely on Internet. We will use it.
\cite{ChobChob} is the basis of dissertation \cite{Ciob3}, that is
why let us analyze \cite{ChobChob} in more details.

$10a_1$. \cite{ChobChob} investigates the algebras with continuous
signature ($E$-algeb\-ras), which are defined in items 1b) -- 1e).
No property of continuous signature is considered in
\cite{ChobChob} . According to item 1e) this notion may be
ignored. In work \cite{ChobChob} the notion of continuous
signature is introduced only to confuse the reader.

10b). In "Introduction" it is mentioned that \textmd{one of the
main problem examined in the present article is the following. Let
$X$ be a subspace of the space $Y$, $K$ be a class of topological
algebras, $F(X,K)$ and $F(Y,K)$ be the free topological algebras
of the spaces $X$ and $Y$, respectively. Assume that $X \subseteq
F(X,K)$, $Y \subseteq F(Y,K)$ and $\varphi: F(X,K) \rightarrow
F(Y,K)$ is the homomorphism for which $\varphi(x) = x$ for any $x
\in  X$. Under which conditions is $\varphi$ an embedding?}

$10b_1)$. According to Theorem 7, this question has a meaning if
and only if the class $K$ is an almost quasivariety.

$10b_2)$. From definition of free algebra given by defining space
(item 10g)) and  $X \subseteq F(X,K)$, $Y \subseteq F(Y,K)$ it
follows that $X$, $Y$ are a $K$-free generators of $F(X,K)$ and
$F(Y,K)$ respectively. Then by Lemma 3 $X$, $Y$ are a
$K$-independent generators. According to "$\varphi(x) = x$ for any
$x \in X$" we consider $X$ as subset of $Y$. Then from the
assertion after Corollary 13 it follows that $\varphi(X)$ is a
$K$-independent set. Again, by Lemma 3 $\varphi(X)$ generates a
$K$-free system. Finally, from item 2) of Theorem 8 a positive
answer to question 10b) follows.

$10b_3)$. The question from item 10b) is almost not considered at
all in the analyzed paper. It is not clear why this question was
stated in the "Introduction".

Now we present literally the excerpts 10c), 10d), 10e), 10f),
10g), 10h).

10c). \textmd{A class $K$ of $E$-algebras is called a quasivariety
of $E$-algebras if the following properties hold:}

\textmd{1Q. a Cartesian product of $E$-algebras from $K$ is an
$E$-algebra from $K$;}

\textmd{2Q. if $A \in  K$ and $B$ is a subalgebra of $A$, then $B
\in  K$.}

$10c_1)$. The definition 10c) is not correct; it differs from the
classical definition of quasivariety of algebraic systems,
presented in Proposition 3.

10d). \textmd{A class $K$ of topological $E$-algebras is called a
topological quasivariety of topological $E$-algebras if the
following properties hold:}

\textmd{1TQ. the topological product of topological algebras from
$K$ is a topological $E$-algebra from $K$;}

\textmd{2TQ. if $A \in K$ and $B$ is a closed subalgebra of $A$,
then $B \in K$.}

$10d_1)$. Recall (\cite{Mal3}) that an algebraic system for which
the basic set of elements is a topological space and the basic
operations are continuous is called topological algebraic system.
Then from here  it follows that the definition 10d) is not
correct, as it contradicts the definition from item 10c). In the
condition of 2TQ, it is necessary to require that $B$ be a
subalgebra, but not a closed subalgebra.

10e). \textmd{A topological quasivariety $K$ of topological
$E$-algebras is a complete topological quasivariety if the
following property holds:}

\textmd{3TQ. If $A \in K$ and a topological $E$-algebra $B$ is
isomorphic with $A$, then $B \in K$.}

$10e_1)$. In the literature,, see, for example, \cite[pag.
207]{Mal1}, \cite{Con}, any class of algebraic systems that
satisfies the condition 3TM is called abstract. By \cite[page
267]{Mal1} any quasivariety or variety of algebraic systems is an
abstract class. Then any class of algebraic systems that satisfies
any of the   conditions $V_1$ -- $V_3$ of Proposition 2 or $Q_1$
-- $Q_4$ of Proposition 3 is abstract. Hence the notion introduced
in item 10e) is without sense. The notion of complete quasivariety
coincides with notion of quasivariety.

10f). \textmd{If the signature $E$ is discrete and $K$ is a
complete topological quasivariety, then $K$ is a variety of
$E$-algebras.}

$10f_1)$. The assertion 10f) is false. is a complete non-sense.
From item $10a_1$ it follows that the discrete or non-discrete
signature does not influence the topological and algebraical
structure of topological quasivariety. Then from items $10e_1)$
and 10f) it follows that any topological quasivariety of
topological $E$-algebras is a variety of $E$-algebras. It is not
clear,  is a variety of topological $E$-algebras or discrete
$E$-algebras?

\textmd{\textbf{Definition 3.1.} Let $K$ be a class of topological
$E$-algebras and $X$ be a space.}

10g). \textmd{(\textbf{3.1T}). A couple $(F(X,K), i_X)$ is called
a topological free algebra of the space $X$ in the class $K$ if
the following conditions hold:}

\textmd{(1) $F(X,K) \in K$ and $i_X: X \rightarrow F(X,K)$ is a
continuous mapping;}

\textmd{(2) the set $i_X(X)$ topologically generates $F(X,K)$;}

\textmd{(3) for each continuous mapping $f: X \rightarrow G \in K$
there exists a continuous homomorphism $\overline f: F(X,K)
\rightarrow G$ such that $f = \overline f \circ i_X$. The
homomorphism $\overline f$ is called the homomorphism generated by
$f$.}

$10g_1)$. By Theorem 7, the topological free algebra of the space
$X$ in the class $K$ exists if and only if the class $K$ is an
almost quasivariety of topological algebras. In such case the
condition (2) if item 10g) has the the form: ($2^{\prime}$) the
set $i_X(X)$ algebraically generates $F(X,K)$.

10h). \textmd{(\textbf{3.1A}). A couple $(F^a(X,K), j_X)$ is
called an algebraically free algebra of the space $X$ in the class
$K$ if the following conditions hold:}

\textmd{(1) $F^a(X,K)$ is a subalgebra of some topological algebra
$G(X) \in K$ and $j_X: X \rightarrow F^a(X,K)$ is a mapping;}

\textmd{(2) the set $j_X(X)$ algebraically generates $F^a(X,K)$;}

\textmd{(3) for each mapping $f: X \rightarrow G \in K$ there
exists a continuous homomorphism $\overline f: F^a(X,K)
\rightarrow G$ such that $f = \overline f \circ j_X$.}

$10h_1)$. It is known (see, for example \cite[pag. 66]{Eng}) that
the notion of mapping of topological spaces is not a topological
notion. One of the main and elementary notion of topological
theory is the notion of continuous mapping. From conditions "$j_X:
X \rightarrow F^a(X,K)$ is a mapping" and "$f: X \rightarrow G \in
K$ is a mapping" it follows that the introduced notion of
algebraically free algebra $F^a(X,K)$ is not a topological notion.
This situation may be corrected only by: a) considering that the
mappings $j_X$, $f$ are continuous; b) considering $j_X$, $f$ as
mappings and that $X$ is a discrete space, as every mapping of
discrete space is continuous \cite[Example 1.4.2]{Eng}. According
to item $10g_1)$ the case a) is the definition \textbf{3.1T}. The
case b) is not possible, as discrete spaces do not meet the
conditions $K_1$, $K_2$ (see the assertions after conditions
$K_1$, $K_2$).

$10h_2)$. Consequently, the introduced notion of algebraically
free algebra from item 10h) is senseless and this notion cannot be
used.

10i). \textmd{In \cite{Chob12}, \cite{Chob3}, \cite{Chob2} the
following theorem was proved.}

$10i_1)$. \textmd{\textbf{Theorem 3.2.} Let $K$ be a topological
quasivariety of topological $E$-algebras and $X$ be a non-empty
space.}

\textmd{1. The free object $(F(X,K), i_X)$ exists and it is
unique.}

\textmd{2. The free object $(F^a(X,K), j_X)$ exists and it is
unique.}

\textmd{3. There exists a unique continuous homomorphism $k_X:
F^a(X,K) \rightarrow F(X,K)$ such that $i_X = k_X \circ j_X$ and
$k_X(F^a(X,K))$ is a dense subalgebra of the algebra $F(X,K)$.}

\textmd{4. If $K$ is a quasivariety, then $F^a(X,K) \in K$,
$k_X(F^a(X,K)) = F(X,K)$ and the set $i_X(X)$ generates the
algebra $F(X,K)$.}

\textmd{5. If $K$ is a non-trivial complete quasivariety (or a
non-trivial quasivariety and an $\omega$-class), then:}

\textmd{- $i_X$ is an embedding of $X$ into $F(X,K)$ and $k_X$ is
a continuous isomorphism of $F^a(X,K)$ onto $F(X,K)$;}

\textmd{ - if $X$ and $E$ are $k_{\omega}$-spaces, then $F(X,K)$
is a $k_{\omega}$-space, too.}

$10i_2)$. Item 10i) contains a misleading statement. The journal
where the work \cite{Chob3} was published is not indicated
correctly. For this see item 1t).

$10i_3)$. The theorem from item $10i_1)$ is not proved at all in
the papers mentioned  in item 10i). For this see items 1j) -- 1l),
and also item 2d).

$10i_4)$. The assertion 1 from item $10i_1)$ is not correct in
such a form. The correct assertion can be found in Theorem 7.

$10i_5)$. The assertions 2 -- 4 from item $10i_1)$ are senseless
according to item $10h_2)$.

$10i_6)$. Assertion 4 and, according to item 10e), assertion 5
consider the quasivarieties, which by definition from item 10c)
contains only (discrete) $E$-algebras. But these assertions
consider a topological free algebra $F(X, K)$, which is not
correct.

10j). Section "2. PRELIMINARIES" introduces the notions of uniform
signature, uniform $E$-al\-gebra, uniform quasivariety of uniform
$E$-algeb\-ras, complete uniform quasivariety, $\omega$-class
similarly of notions of uniform signature, uniform $E$-algebra,
uniform quasivariety of uniform $E$-algebras, complete uniform
quasivariety, $\omega$-class respectively (only the word
"continuous" is replaced by "uniform").

$10j_1)$. The class of uniformities spaces coincides with the
class of completely regular spaces \cite[Theorem 8.1.20]{Eng}.
Then from item $10a_1$ (see, also, item 1e)) it follows that the
notion of uniform signature should be ignored.

$10j_2)$. The introduced notions from item $10j_1)$ are
investigates in sections 7 and 8. The content of this section is
not Mathematics. To confirm the above-mentioned we will show some
excerpts.

$10j_3)$  (pag. 67). \textmd{If the signature $E$ is discrete and
$K$ is a complete quasivariety of weakly uniform topological
$E$-algebras, then $K$ is a weakly uniform quasivariety.}

$10j_4)$. (Proof of Theorem 7.8). \textmd{Since $K$ is a complete
quasivariety, the homomorphism $k_X: F^a(X,K) \rightarrow F(X,K)$
is a continuous isomorphism for any space $X$. Moreover, if $K'$
is the variety generated by the class $K$, then $F(X,K) = F(X,K')$
for any space $X$ \cite{Chob12}, \cite{Chob3}, \cite{Chob2}. Thus
we can suppose that $K$ is a variety, i.e. there exists a family
of identities $\mathcal{J}$ such that $K$ is the class of
topological $E$-algebras with the identities. In this case $(cG,
\mathcal{T}(cW_G)) \in K$ for any $G \in K$ \ldots The assertion 1
is proved. The assertion 2 follows from the assertion 1. The
assertion 3 follows from \cite{ChobKir4}. The proof is complete.}

$10j_5)$. This excerpt even contains a link to an non-existing
paper $10j_2)$ (see also item $10i_2)$). In the first sentence of
item $10j_4)$ the senseless assertion from item 10f) is used, as
well as the senseless definition of $F^a(X,K)$ from item 10h) (see
and the item $10h_2)$). There is no need to further comment item
$10j_2)$.\vspace*{0.1cm}

10k). The reasoning from item $10j_2)$ is also valid for the other
sections of the analyzed paper \cite{ChobChob}. Let us make a
short analysis.\vspace*{0.1cm}

10l).  In section 4 the topological free algebra $F(X, K)$ is
investigated. Then from Theorem 7  we should
 the class $K$ to almost quasivariety  in item
9g). In such a case $K$ contains free topological algebraic
systems $F_m$ of any given rank $m \geq 1$ by Theorem 7.

$10l_1)$.  We consider that $K$ is an almost quasivariety of
topological algebraic systems given by defining space and defining
relations. Then from Definition 1 it follows that the condition
"\ldots and $\alpha_G: G \rightarrow \alpha\beta_{\mathcal{P}}F(X,
K)$ is an embedding \ldots" from Theorem 4.8 takes the shape of
topological equality: $G = \alpha\beta_{\mathcal{P}}F(X,  K)$.

$10l_2)$.  From item 2) of Theorem 8  and definition of \break
$\alpha\beta_{\mathcal{P}}F(X, K)$ the topological equality $F(X,
K) = \alpha\beta_{\mathcal{P}}F(X, K)$ follows.

$10l_3)$.  From item 8a) and definition of $(\beta_{\mathcal{P}},
\beta_X)$ it follows that $(\beta_{\mathcal{P}}, \beta_X)$ is the
free topology regarding to space $X$.

$10l_4)$.  The proof of Theorem 4.8 is senseless. Without going
into details, compare the items $10l_2)$, $10l_3)$ and the
expression from Proof: "The assertion 2 is proved. The assertion 3
follows from the assertion 2".

$10l_5)$. The conditions of Corollary 4.9 is completely
senseless.\vspace*{0.1cm}

10m). Now we present literally the excerpts $10m_1)$ -- $10m_5)$
from Section 5.

$10m_1)$. \textmd{\textbf{Corollary 5.5.} If $X$ is a discrete
space, then $F(X,K)$ is a projective algebra in the class $K$.}

$10m_2)$. \textmd{\textbf{Theorem 5.7.} For a topological
$E$-algebra $G \in K$ the following assertions are equivalent:}

\textmd{1. $G$ is a projective algebra in the class $K$;}

\textmd{2. if $X$ is a discrete space and $|G| \leq |X|$, then $G$
is an $\alpha$-retract of the algebra $F(X,K)$;}

\textmd{3. there exists a discrete space $X$ such that $G$ is an
$\alpha$-retract of the algebra $F(X,K)$.}

$10m_3)$. \textmd{\textbf{ Proof} \ldots  Since the class $K$ is
not trivial and $X$ is discrete, $i_X: X \rightarrow F(X,K)$ is an
embedding \ldots}

$10m_4)$. \textmd{ A topological quasivariety $K$ is a Schreier
class if for every discrete space $X$ the subalgebras of $F(X,K)$
are free in $K$.}

\textmd{From Theorem 5.7 it follows}

$10m_5)$. \textmd{\textbf{Corollary 5.9.} If $K$ is a Schreier
class, then only the free algebras $F(X,K)$ of discrete spaces $X$
are projective in $K$.}

$10m_6)$. Items $10m_1)$ -- $10m_5)$ consider a topological free
$E$-algebra \break  $F(X, K)$ of quasivariety of topological
$E$-algebras with discrete space $X$. But according to item
$10h_2)$ such topological free $E$-algebra $F(X, K)$  do not
exist. Hence, the assertions of items $10m_1)$ -- $10m_5)$, as
well as all results from Section 5, are false and senseless.
Moreover, the absurdity of  item $10m_3)$ will be discussed
bellow.

$10m_7)$. Let us show that a topological free $E$-algebra $F(X,
K)$ of quasivariety of topological $E$-algebras with discrete
space $X$ does not exist. Really, as in item 9g) let us transform
the class $K$ in almost quasivariety in the sense of Lemma 1.
According to Theorem 7  only and only the almost quasivariety
contain topological free algebras $F(X, K)$. Then from Lemma 1 it
follows that the space $X$ should meet the conditions $K_1)$,
$K_2)$ from Section 2. After conditions $K_1)$, $K_2)$ it is shown
that the discrete spaces do not meet the conditions $K_1)$,
$K_2)$, they only generate such spaces.

Now we pass to Section 6, we quote.

10n). \textmd{Let $E$ be a continuous signature.}

\textmd{Fix a non-trivial quasivariety $K$ of topological
$E$-algebras. Suppose that $i_X: X \rightarrow F(X,K)$ is an
embedding}

$10n_1)$. \textmd{and $k_X: F^a(X,K) \rightarrow F(X,K)$ is a
continuous isomorphism for every space $X$. In this case we can
consider that $X = i_X(X)$ is a subspace of $F(X,K)$ and $a(X,
F(X,K)) = F(X,K)$ for any space $X$.}

$10n_2)$. \textmd{\textbf{Corollary 6.2.} $X$ is a closed subspace
of the space $F(X,K)$ for any space $X$.}

$10n_3)$. The assertion from item $10n_2)$ is false. The space $X$
should meet the conditions $1TQ$, $2TQ$ from item 10d).

 $10n_4)$. According to item $10h_2)$ the assumption from
item $10n_1)$ senseless. Therefore the Proposition 6.1 and the
Corollary 6.2 from item $10n_3)$ are also senseless.

10o). The remaining part of Section 6, Proposition 6.3 and
Corollaries 6.3, 6.4 consider the case when the signature $E$ is
discrete and $K$ is a complete quasivariety. By item 10f) $K$ is a
variety in sense of Academician Choban. It is not clear where the
properties of variety are used. But they use the senseless
Proposition 6.1 and Corollary 6.2. Hence, this part of Section 6
is lacking any sense.

10p). On page 76 it is written that in Examples 12.1 -- 12.8 some
important classes $K$ of universal algebras are mentioned. These
examples are senseless. It is not even the case to comment them.
We mention only that these examples contain the statement: the
mapping $i_X: X \rightarrow F(X, K)$ is an embedding for any space
$X$ ($X$ is the defining space for topological free algebra $F(X,
K)$ of class $K$ according to item 10g)). The last statement has a
direct relation to question B) from section 2. Comments and the
sufficiency condition for a positive answer to this question can
be found in \cite[Remark 2]{Mal3}.

10q). For a short analysis of the remaining sections 9 -- 13 we
quote the item $10q_1)$ from \cite{Mal2} and the items $10q_2)$
and $10q_3)$ from \cite{RezUsp}.

$10q_1)$. \textmd{A class $K$ of topological $E$-algebras is a
Mal'cev class if there exists a ternary derivate operation $m$
such that $m(x, x, y) = m(y, x, x) = y$ for any $\mathcal{A} \in
K$ and any $x, y \in \mathcal{A}$. In this case we can consider
that  $m \in E_3$ and we say that $m$ is a Mal'cev ternary
operation.}

$10q_2)$. \textmd{A Mal'cev operation on a space $X$ is a
continuous function  $f: X^3 \rightarrow X$ such that $f(x, y, y)
= f(y, y, x) = x$, for all $x, y \in X$. A topological space is
called  Mal'cev if it admits a Mal'cev operation.}

$10q_3)$. (Theorem 1.6). \textmd{Let $X$ be a pseudocompact
Mal'cev space. Then:}

\textmd{(a) the Stone-Cech compactification $\beta X$ is
Dugundji;}

\textmd{(b) every Mal'cev operation on $X$ extends to a Mal'cev
operation on $\beta X$}

\textmd{(c) $X$ is a retract of a topological group.}

$10q_4)$. (\cite{ComRos}). \textmd{Let $G$ be a pseudocompact
topological group. Then the group operation on $G$ extends to a
continuous group operation on $\beta G$ which makes $\beta G$ into
a topological group.}

$10q_5)$. The definition and meaning of operation $m$ from item
$10q_1)$ is not clear. On the one side, $m$ is a ternary derivate
operation, and on the other side $m \in E_3$ is a basic operation
of topological $E$-algebra. A correct definition and meaning of
Mal'cev operation are shown in Theorems 1, 3. The significant
difference between the definitions from items $10q_1)$ and
$10q_2)$ are described in details in item $^{(8)}$.

$10q_6)$. In sections 9 -- 13 the authors are used a statement for
topological $E$-algebras with nonempty set of basic operation,
which if not wrong, than it least is not proved: they transfer
theorem from item $10q_3)$ for topological spaces on topological
$E$-algebras. This proves the assertions from items $10q_3)$ and
$10q_4)$.

$10q_7)$. The unproved statement from item $10q_6)$ has the
following form. Let $G$ be a pseudocompact topological
$E$-algebra. Then the basic  operations  on $G$ extend to a
continuous basic operations on $\beta G$ which makes $\beta G$
into a topological $E$-algebra.

$10r)$. As proved partially above, we confirm again that sections
9 -- 13 do not contain any mathematical assertion proven
correctly. If Mal'cev algebra is considered in sense of item
$10q_2)$,   then well known assertions are listed (Theorem 9.1,
Theorem 12.9, but item 3 is not clear.). In some other cases
sections 9 -- 13 contain senseless assertions. We confirm this
only partially and quote literally Proposition 9.2 and Theorem
9.3.

$10r_1)$. \textmd{\textbf{Proposition 9.2.} Let $G$ be a
pseudocompact Mal'cev $E$-algebra. If $G$ is a dense subalgebra of
a Mal'cev $E$-algebra $A$, then $G \subseteq A \subseteq \beta
G$.}

$10r_2)$. \textmd{\textbf{ Proof.} On $\beta A$ and $\beta G$
there exist the structures of Mal'cev algebras such that:}

\textmd{- $A$ is a subalgebra of the algebra $\beta A$;}

\textmd{- $G$ is a subalgebra of the algebra $\beta G$;}

$10r_3)$. \textmd{- there exists a continuous homomorphism $g:
\beta G \rightarrow \beta A$ such that $g(x) = x$ for any $x \in
G$.}

\textmd{The homomorphism $g$ as a quotient mapping of the compact
Mal'cev algebra $\beta G$ onto the Mal'cev algebra $\beta A$ is an
open mapping. Thus $g$ is an isomorphism and $\beta A = \beta G$.
The proof is complete.}

$10r_4)$. \textmd{\textbf{Theorem 9.3.} Let $K$ be a topological
quasivariety of compact $E$-algebras, $K$ be a Mal'cev class, $X$
be a subspace of the free $E$-algebra $F(X, K)$ and $B$ be a
pseudocompact subalgebra of $F(X, K)$ such that $X \subseteq B$.
Then $\beta B = F(X, K)$.}

$10r_5)$. \textmd{\textbf{Proof.} By virtue of Proposition 9.2 we
have $B \subseteq  F(X, K)$ and $\beta B = \beta F(X, K) = F(X,
K)$. The proof is complete.}

$10r_6)$. From items $10q_6)$ and $10q_7)$  it follows that the
algebras $\beta A$ and $\beta G$ of item $10r_2)$ are Mal'cev in
sense $10q_2)$ by item $10q_4)$.

$10r_7)$. We will not analyze the topological side of $10r_3)$.
Without this we conclude that Proposition 9.2 is proved for
Mal'cev algebras in sense of item $10q_2)$.

$10r_8)$. The equalities $\beta B = \beta F(X, K) = F(X, K)$ from
item $10r_7)$ are false. By item $10r_7)$ $\beta B$ is a Mal'cev
algebra in sense of item $10q_2)$, but $F(X, K)$ is a Mal'cev
algebra in sense of item $10q_1)$.

$10r_9)$. According to item $10r_8)$ the Theorem 9.3 is false. The
Proposition 9.2 and the Theorem 9.3 are a mechanical
transformation of Pestov's result mentioned in \cite[pag.
87]{RezUsp}. The result is senseless.

$10s)$. Similarly to item $10r_9)$ for the results of Section 9 it
is possible to ascertain also all results of Sections 10 -- 13. We
present only three quote for confirmation.

$10s_1)$. \textmd{\textbf{Theorem 11.2.} Let $G$ be a
pseudocompact Mal'cev $E$-algebra \ldots}

\textmd{\textbf{Proof.} From Reznicenko-Uspenskiy theorem
\cite{RezUsp} it follows that $G$ is a subalgebra of the compact
Mal'cev $E$-algebra $\beta G$ \ldots} According to items items
$10q_6)$ and $10q_7)$ $\beta G$ is not  a Mal'cev algebra in sense
of item $10q_1)$.

$10s_2)$. \textmd{\textbf{Corollary 11.3.} Let $K$ be a
topological variety of compact $E$-algebras \ldots} About what
reliability of Corollary 11.3 there can be a speech  if by item
10f) the authors  not own elementary notion  notion of variety of
algebras.

$10s_3)$. Sections 9 -- 13 contain referrals to paper
\cite{Chob2}. We showed above that the paper \cite{Chob2} consists
only of false assertions, that it is senseless.

$10t)$.  We confirm again that according to the above-mentioned,
paper \cite{ChobChob} is a non-sense presented on 30 journal
pages. \vspace*{0.1cm}

11a). Paper \cite{Chob7} is written under the influence and in the
manner of Academician Choban M. M. works. The main result (Theorem
1.10) is not correct, as false notions and results of Choban M. M.
are used.

$11a_1)$. First, it is necessary to ignore the notions of
continuous and discrete signature according to item 1e). The
notion of Tychonoff product (Definition 1.3) is not clear. For the
correct notion of Tychonoff (or direct) product of topological
algebras see, for example, \cite[Proposition 2.3.1]{Eng} (or
\cite[pag. 174]{Mal3}). Further, we quote the items $11a_2)$ --
$11a_3)$.

$11a_2)$. \textmd{\textbf{Definition 1.5.} A class $\mathcal{K}$
of topological $E$-algebras is called a \textit{complete
quasivariety} if}

\textmd{(1) for each topological $E$-algebra $X \in \mathcal{K}$,
each $E$-subalgebra of $X$ belongs to the class $\mathcal{K}$;}

\textmd{(2) for any topological $E$-algebras $X_{\alpha} \in
\mathcal{K}$, $\alpha \in  A$, their Tychonov product
$\prod_{\alpha \in A} X_{\alpha}$  belongs to the class
$\mathcal{K}$;}

\textmd{(3) a Tychonov $E$-algebra belongs to $\mathcal{K}$ if it
is algebraically isomorphic to a topological $E$-algebra $Y \in
\mathcal{K}$.}

$11a_3)$. The definition of item $11a_2)$ is false, such a
nonempty classes of  complete quasivarieties does not exist.
Really, assume that by item $11a_1)$ the notion of Tychonoff
product is defined correctly. Similarly to items 9d) -- 9g) we
transform by Lemma 1 the class $\mathcal{K}$ from item $11a_2)$
into an almost quasivariety of topological algebras with given
defining space $X$ with condition (3). The space $X$ should
satisfy the conditions $K_1)$, $K_2)$ from Section 2. But the
Tychonoff spaces from condition (3) do not satisfy the conditions
$K_1)$, $K_2)$. Hence the condition (3) it without sense.

$11b_1)$. \textmd{\textbf{Definition 1.6.} Let $\mathcal{K}$ be a
complete quasivariety of topological $E$-algebras. A \textit{free
topological $E$-algebra in $\mathcal{K}$ over a topological space
$X$} is a pair $(F_{\mathcal{K}}(X), \eta)$ consisting of a
topological $E$-algebra $F_{\mathcal{K}}(X) \in \mathcal{K}$ and a
continuous map $\eta: X \rightarrow F_{\mathcal{K}}(X)$ such that
for any continuous map $f: X \rightarrow Y$ to a topological
$E$-algebra $Y \in \mathcal{K}(X)$ there is a unique continuous
$E$-homomorphism $h: F_{\mathcal{K}}(X) \rightarrow Y$ such that
$f = h \circ \eta$.}

$11b_2)$. The definition of free topological $E$-algebra from item
$11b_1)$ is not complete, not coinciding with corresponding
definitions from \cite{Mal3}, \cite{Chob3}, \cite{DC},
\cite{Chob9}, \cite{Chob19}, \cite{ChobChob}: see Definition 1
from Section 2, items 1h), 2b), $4d_3)$, 6d). In Definition 1.6 it
is necessary to add the condition: the set $\eta(X)$ topologically
(algebraically) generate $(F_{\mathcal{K}}(X), \eta)$. But then we
obtain a contradiction with the non-senses from Theorem 1.7.

 $11c_1)$. \textmd{The construction $(F_{\mathcal{K}}(X)$ of a free topological
$E$-algebra has been intensively studied by M.M.Choban
\cite{Chob9}, \cite{Chob19}.}

$11c_2)$. Item  $11c_1)$ contains the most severe error of the
analyzed paper \cite{RezUsp}. We showed earlier (items $4e_4)$,
$6e)$) that the Choban's works \cite{Chob9}, \cite{Chob19} consist
only of erroneous,  senseless  statements and introduced notions,
these works are  anti-scientific.

 $11c_3)$. \textmd{In particular, he proved that for each complete
 quasivariety $\mathcal{K}$K of
topological $E$-algebras and any topological space $X$ a free
topological $E$-algebra $(F_{\mathcal{K}}(X), \eta)$ exists and is
unique up to a topological isomorphism.}

 $11c_4)$. In works \cite{Chob9}, \cite{Chob19} the statement from item
$11c_3)$ is not proven, see the analyze of these papers.

 $11d_1)$. \textmd{Also he proved the following important result,
 see \cite[2.4]{Chob9}:}

$11d_2)$. \textmd{\textbf{Theorem 1.7} (Choban). If $\mathcal{K}$
 is a non-trivial complete quasivariety
of topological $E$-algebras, then for each Tychonov space $X$ the
canonical map $\eta: X \rightarrow  F_{\mathcal{K}}(X)$ is a
topological embedding and $F_{\mathcal{K}}(X)$ coincides with the
subalgebra $<\eta(X)>$ generated by the image $\eta(X)$ of $X$ in
$F(X, \mathcal{K})$.}

$11d_3)$. The Theorem 1.7 from item $11d_2)$ is false, is a
complete non-sense. In item $11d_1)$ it is mentioned that this is
Theorem 2.4 from \cite{Chob9}. In \cite{Chob9} it is mentioned
that it is proved in \cite{Chob3} and \cite{DC}. According to
items 1l), 2d) such a statement in  \cite{Chob3} and \cite{DC} is
not proved , see the analysis of these papers.

$11d_4)$. \textmd{Since $\eta: X \rightarrow (F_{\mathcal{K}}(X),
\eta)$ is a topological embedding, we can identify a Tychonov
space $X$ with its image $\eta(X)$  in $(F_{\mathcal{K}}(X),
\eta)$) and say that the free $E$-algebra $(F_{\mathcal{K}}(X),
\eta)$) is algebraically generated by $X$.}

$11d_5)$. The corollary of Theorem 1.7 from item $11d_4)$ is
completely false and is essentially used to define the function
$F_{\mathcal{K}}$.

$11e)$. The Theorem 1.8 is the Theorem 4.1.2 from \cite{Chob19}
and is a complete non-sense as Theorem 1.7.

$11f)$. The complete set of Choban's non-senses from items
$11a_1)$, $11a_2)$, $11b_1)$, $11c_3)$, $11d_2)$, $11d_4)$,
$11d_5)$ are used essentially to prove the main result of Theorem
1.10. As corollary the Theorem 1.10 and its proof is a complete
non-sense.

11g). The correct variant of Theorem 1.10 would be the following.
It is necessary to analyze an almost quasivariety of topological
algebras of given signature instead of a complete quasivariety of
topological $E$-algebras of contable discrete signature $E$
$\mathcal{K}$ from Theorem 1.10 according to Theorem 7 from
present paper  and describe the submetrizable (compactly
finite-dimensional $ANR(k_{\omega})$ spaces that satisfy the
conditions $K_1)$ and $K_2)$.\vspace*{0.1cm}

12a).  Paper \cite{ChobChir44}, like paper \cite{ChobChob}, is
published in ROMAI Journal and  like \cite{ChobChob} have the same
scientific importance (see the item $10t$): it is a non-sense,
presented on 25 pages.

$12a_1)$. Paper \cite{ChobChir44} investigates the topological
algebras with continuous signature, in short $E$-algebras. For the
definition of $E$-algebra see item 1b). According to item 1e) the
notion of continuous signature should be ignored.

$12a_2)$. In the beginning of this paper it is mentioned that the
paper discusses some old and new results and open problems. We
will show below that the even the well known results are presented
incorrectly. Moreover, what open problems they are referring to,
if the authors are not aware of the basic notions: polynomial,
term, identity, variety, quasivariety, topological algebra, free
topological algebra and others. To support the above mentioned, we
quote the items $12b)$, $12c)$, $12d)$, $12e)$, $12g)$.

$12b)$. \textmd{Let $G$ be a topological $E$-algebra.}

\textmd{The polynomials are constructed in the following way:}

\textmd{- $E$ are polynomials;}

\textmd{- if $n \in N$, $n \geq 1$, $u \in E_n$, $p_i$ is an
$m_i$-ary polynomial, then $p = u(p_11, ..., p_n)$ is an $m$-ary
polynomial, where} \textmd{$m = m_11 + m_22 + ... + m_n \quad
\text{and}$} \textmd{$p(x_1, ..., x_m) = u(p_1(x_1, ..., x_{m1}),
..., p_n(x_{m_{n-1}+1}, ..., x_m)).$}

$12b_1)$. The definition of polynomial from item $12b)$ is false.
$E = \cup E_n$ is the set of symbols of $n$-ary basic operations
and if $n > 1$ $E_n$ cannot be a polynomial. Correct: see, for
example, \cite{Mal2}, \cite{Mal3} or "Preliminaries. Topological
algebraic systems, quasivarieties" Section.

$12c)$. \textmd{Let $n \geq m \geq 1$, $p$ be an $n$-ary
polynomial and $q : \{1, 2, ..., n\} \rightarrow \{1, ...,m\}$ be
a mapping. Then $v(x_1, ..., x_m) = p(x_{q(1)}, x_{q(2)}, ...,
x_{q(n)})$ is an $m$-ary term. The polynomials are terms too. If
$u$ is an $n$-ary term and $v$ is an $m$-ary term, then $u(x_1,
..., x_n) = v(y_1, ..., y_m)$ is an identity on $E$-algebras.}

$12c_1)$. The definition of term from item $12c)$ is completely
senseless. Correct: see, for example, \cite[pag. 141]{Mal1}  or
"Preliminaries. Topological algebraic systems, quasivarieties"
Section.

$12c_2)$. The definition of identity on $E$-algebras from item
$12c)$ is completely senseless. Correct: see, for example,
\cite[pag. 189, 268, 276]{Mal1} or Proposition 5 from Section 2.

$12d)$. \textmd{A class $\mathcal{K}$ of topological $E$-algebra
is called a $T_i$-quasivariety if:}

\textmd{- any algebra $G \in \mathcal{K}$ is a $T_i$-space,}

\textmd{- if $G \in \mathcal{K}$ and $B$ is a subalgebra of $G$,
then $B \in \mathcal{K}$,}

\textmd{- the topological product of algebras from $\mathcal{K}$
is a topological algebra from $\mathcal{K}$,}

\textmd{- if $(G, \mathcal{T}) \in \mathcal{K}$, $\mathcal{T}'$ is
a $T_i$-topology on $G$ and $(G, \mathcal{T}')$ is a topological
$E$-algebra, then $(G, \mathcal{T}') \in \mathcal{K}$.}

$12d_1)$. An algebraic system for which the basic set of elements
is a topological space and the basic operations are continuous is
called topological algebraic system (\cite{Mal3}). Every mapping
of discrete space is continuous \cite[Example 1.4.2]{Eng}.

$12d_2)$. According to item $12d_1)$ we can consider that the
topology $\mathcal{T}'$ from the last condition of item $12d)$ is
the discrete topology. In such a case the definition of
$T_i$-quasivariety from item $12d)$ coincides with the definition
of topological $Q$-quasivariety from item $9b)$. By item $9i)$ the
definition of topological $Q$-quasivariety from item $9b)$ is a
non-sense. Hence, we conclude.

$12d_3)$. The definition of $T_i$-quasivariety from item $12d)$ is
false and senseless.

$12e)$. \textmd{If $O$ is a set of identities and $V(E, \Omega,
i)$ is the class of all topological $E$-algebras with identities
$\Phi$, which are $T_i$-spaces, then $V(E, \Omega, i)$ is a
\textit{$T_i$-variety}. Any $T_i$-variety is a
$T_i$-quasivariety.}

$12e_1)$. How can the definition of $T_i$-variety from item $12e)$
can be accurate if according to item $12c_2)$ the authors do not
know the definition of identity. This confirms also the
proposition from item $12e)$.

$12f)$. The main topological structure in all assertions from
Sections 4, 5, 7, 9  are  the $T_i$-quasivarieties and their
topological free algebras. Then from items $12d_3)$ and $12e_1)$
it follows that all assertions from Section 4 are false. According
to items  9d) - 9g) we transform by Lemma 1 the
$T_i$-quasivarieties from item $12d)$ into an almost
quasivarieties of topological algebras with given defining space
$X$. In such a case, by Theorem 7 there exist topological free
algebras. Even with such transformations the proofs of theorems
are false. The proofs use free topological algebras with defining
discrete space $X_d$. By item $12d_2)$ this is not allowed.

$12g)$. \textmd{The necessity and sense of defined notions is not
clear: $S$ -semito\-pological $E$-algebra, separately continuous,
$P$-paratopological $E$-algebra, \break ($P, S )$-quasitopological
$E$-algebra (pag. 9), left (right) topological group (pag. 11),
$S$-simple $E$-algebra (pag. 12), topological $T$-groupoid (pag.
14), a classes $V(E, \varphi)$, $V(E, u\varphi)$, $V(E, \Pi)$,
primitive solution (pag. 17), topological bigroupoid,(pag. 18),
topological $E$-automaton (pag. 19).}

$12h)$.  All assertions from the analyzed work, related with the
introduced notions from item $12g)$ are false or trivial. A
sufficient condition for Problem 5.2 is \cite[Proposition
IX.1.6]{Smith1}. The author of this paper is aware of the
researches into $T$-quasigroups (linear quasigroups). However,
Section 5 contains only some unsuccessful definitions of
$T$-groupoids. Also compare with \cite[Proposition
IX.1.2]{Smith1}.

$12i)$. Many false results: Theorems 4.1 -- 4.3, 7.2 and others
are quoted from other works: \cite{Chir1}, \cite{Chob2},
 \cite{Chob191}, \cite{Chob28}, \cite{ChobKir39},
\cite{ChobKir40}, \cite{ChobChir36}. Paper \cite{ChobChir36} is
not published in the indicated journal.

 $12j)$. Section 2 contains a review of the researches into Mal'ces
 algebras and related notions of homogeneous algebras, of biternary
 Mal'cev algebras, of rectifications, of retracts belonging mainly to Choban M.
M.: \cite{Chob1}, \cite{Chob2}, \cite{AC1}, \cite{AC2},
\cite{ChobChob}, as well as \cite{Chob28},  \cite{Usp1},
\cite{Usp2}, \cite{RezUsp}.

$12j_1)$. All assertion regarding Mal'cev algebras, mentioned in
item $12j)$, are false. These details are indicated when analyzing
works \cite{Chob1}, \cite{Chob2}, \cite{AC1}, \cite{AC2},
\cite{ChobChob}. The most severe mistake: the definition of
Mal'cev algebra does not coincide with notion of Mal'cev algebra
from \cite{Mal2}. This is shown in items $^{(8)}$, $10q_1)$,
$10q_2)$, $10q_5)$.\vspace*{0.1cm}

$13a)$. The notion of free universal algebra ($\Omega$-algebra) is
classical. Practically any textbook about algebra contains a
theory of such notion, see, for example, \cite{Mal1},
\cite{Kurosh}, \cite{Con}, \cite{Smith1}.

$13a_1)$. Let us consider an universal algebra with operations
$\Omega$ and let $X = \{x, y, \ldots\}$ to be an auxiliary
nonempty set, whose elements will be called free elements. Let us
assign some symbols to all nullary operations from $\Omega$, if
any, which will be symbols of nullary operations. We denote them
by $0_1, \ldots, 0_r$.

$13a_2)$. Let us define the notion of word. Words are all the free
elements and all symbols of nullary operations. If the expressions
$w_1, w_2, \ldots, w_n$ are words, then for every $n$=ary
operation $\omega \in \Omega$, where $n \geq 1$, the formal
expression $\omega(w_1, w_2, \ldots, w_n)$ is also regarded as a
word. The set of all possible words with respect to set of
operations $\Omega$ and the set of all free elements $X$ form the
algebra of words $S(\Omega, X)$.

$13a_3)$. We will say that in universal algebra $G$ with set of
operations $\Omega$ an identity (identical relation) $w_1 = w_2$
is realized, if this equality holds in $G$ when replacing the free
elements from $w_1, w_2$ by arbitrary elements from $G$. The
symbols of nullary operations $0_1, \ldots, 0_r$ are replaced by
these elements from $G$, which mean in $G$ the nullary operations.

$13a_4)$. Let $\Lambda$ be a set of identical relations of form
$w_1 = w_2$ from item $13a_3)$.  The words $v^{\prime},
v^{\prime\prime} \in S(\Omega, X)$ are called equivalent with
respect to $\Lambda$  if one can be transformed into the other by
a finite sequence of transformations of form: let $\Lambda$
contain the identical relation $w_1 = w_2$; let as replace the
free elements contained into it $x_i$, $i = 1, 2, \ldots, k$, by
some words, after which the left and right parts of $w_1 = w_2$
are transformed in words $\overline{w_1}$ and $\overline{w_2}$; if
$\overline{w_1}$ (or $\overline{w_2}$) is a subword of word
$v^{\prime}$, then it is replaced in  $v^{\prime}$ on
$\overline{w_2}$ (or respectively on $\overline{w_1}$).

$13a_5)$. The equivalence from item $13a_3)$ is a congruence
$\sim$ of algebra of words $S(\Omega, X)$. The quotient algebra
$S(\Omega, X)/\sim$ is a free $\Omega$-algebra with free
generators $X$ of class of $\Omega$-algebras with defining
relations $\Lambda$. In effect, for $S(\Omega, X)/\sim$ the set of
free generators is not $X$, but the set of classes of congruence
$\sim$.

$13a_6)$. When performing calculations in free groups, free loops,
free commutative quasigroups and loops, free $TS$-quasigroups and
loops, free Steiner quasigroups, free commutative $IP$-loops and
many other related classes, no harm results if we treat the words
as the actual elements rather than representatives of congruence
classes. Such calculations are greatly facilitated by the
following ideas. We define a word to be in normal form or to be a
reduced word if no reductions of it is possible. In fact, we can
say much more. Any word has a unique reduced form and two words
are congruent if and only if they have the same reduced form
\cite{Smith1}.

$13b)$. Now we pass to the analysis of monograph \cite{Chir1}. We
shall note that cite{Chir1} literally is party of dissertation
\cite{Chir2}, the numbering of chapters  differs only.

$13b_1)$. Similarly, \cite{Chob7} the monograph \cite{Chir1} is
written under influence and in the manner of Academician Choban M.
M. As it uses false notions and results of Choban M. M., similarly
with \cite{Chob7}, all introduced notions and the stated results
from \cite{Chir1} are not correct, are without sense. For the
introduced notions of term, polynomial, identity see the items
12b), 12c). Moreover, it introduced the senseless: let
$\omega_1(x_1, \ldots, x_n) = \omega_2(y_1, \ldots, y_m)$ be an
identity. If $y_1 \notin \{x_1, \ldots, x_n\}$, then $y_1$ is
called a free variable of the identity $\omega_1 = \omega_2$.

$13b_2)$. Monograph \cite{Chir1} introduces and investigates a
senseless definition that in \cite{ChobChir44}: see the items
$12a_1)$, 12b), 12c). Now we present literally the items $13c)$,
$13c_1)$.

$13c)$. \textmd{\textbf{Definition 1.2.6.} A class $K$ of
topological $E$-algebras which are also $T_i$-spaces is called a
complete $T_i$-variety if the following conditions are fulfilled:}

\textmd{(M1). $K$ is closed with respect to Tychonoff products.}

\textmd{(M2). $K$ is closed with respect to subalgebras.}

\textmd{(M3). If $(G; \tau) \in K$ and $(G; \tau^{\prime})$ is a
topological $E$-algebra and also a $T_i$-space, then $(G;
\tau^{\prime}) \in K$.}

(M4). If ($(G; \tau) \in K$ and $(G^{\prime}; \tau^{\prime})$ is a
topological $E$-algebra, $T_i$-space and there exists a continuous
homomorphism $f: G \rightarrow  G^{\prime}$, then  $(G^{\prime};
\tau^{\prime}) \in K$.

$13c_1)$. \textmd{A class $K$ of topological $E$-algebras is
called: $T_i$-quasivariety if conditions M1, M2 hold;
$T_i$-variety if conditions M1, M2, M4 hold; complete
$T_i$-quasivariety if conditions M1 - M3 hold.}

$13c_2)$. According to item 9i) the introduced notions from 13c),
$13c_1)$ are senseless.

$13c_3)$. Similarly to items 10g), 10h) for $T_i$-quasivarieties
$K$ the notions of topological free algebra $F(X, K)$ and
algebraically free algebra $F^a(X, K)$ are introduced, the false
assertions from items 10i), $10i_1)$ -- $10i_6)$ are presented.

$13c_4)$. The false definitions and assertions from items
$13b_1)$, $13c)$, $13c_1)$, $13c_3)$ are basic for more than half
of the analyzed monographs.

$13d)$. The above-mentioned notions and results are part of
Chapter 1. We pass to Chapter 2. Literally we quote the items
$13d_1)$ -- $13d_4)$.

$13d_1)$. \textmd{\textbf{Definition 2.1.1.} Let $E$ be a
continuous signature. The class $K$ of $E$-algebras is called a
Mal'cev class if there exists a polynomial $p(x,y,z)$ such that
the equations $p(x,y,y) = p(y,y,x) = x$ hold identically in $K$.}

$13d_2)$. \textmd{Let $K$ be a Mal'cev class and a complete
$T_i$-variety, $i \geq 0$ … Let $X$ be a Tychonoff space and $F(X,
K)$ be a free topological $E$-algebra of the space $X$ in a given
class $K$.}

$13d_3)$. \textmd{We consider that $X \subset F(X, K) = G$.
Furthermore, we set} \break \textmd{$F_0((X, K)  = X$, $F_1(X, K)
= \sum\{e_{iG}(E_i \times (F_0((X, K))^i : i \leq 3\}$ and for any
$m \geq 2$,
$$ F_m(X, K) = \cup\{e_{iG}(E_i \times (F_{m-1}(X, K))^i : i \leq
m+2\}.$$} $13d_4)$. \textmd{It is clear that $F_m(X, K) \subset
F_{m+1}(X, K)$ for every $m \in N$ and $F(X, K) = \cup\{F_m(X, K)
: m \in N\}$.}

$13d_5)$. From items $12b)$. it follows that $p \in E_3$ for
polynomial $p$ of item $13d_1)$. In such a case there is not
relation between ternary operation $p$ and other basic operations
from $E$. For this see item $10q_5)$. Further, instead of "…the
equations $p(x,y,y) = p(y,y,x) = x$..."  one should read "…the
equalities $p(x,y,y) = p(y,y,x) = x$..."

$13d_6)$. According to item $13c_2)$ the assertion from item
$13d_2)$ is false. Complete $T_i$-varieties do not exist.

$13d_7)$. Item $13d_4)$ is unclear:  how come $F(X, K)$ is a free
algebra. For that see item $13a_4)$. We suppose that $F(X, K)$ is
a free algebra. Then from $F(X, K) = \cup\{F_m(X, K) : m \in N\}$
it follows that every element of $F(X, K)$ has a reduced form.
Then this again contradicts the item $13a_4)$.

$13d_8)$. All proofs of Chapter 2 are based on assertion of form.
Let \break $x_1, x_2, \ldots , x_{2n+1} \in X$. Then $p_n(x_1,
x_2, \ldots , x_{2n+1} \notin F_m(X, K)$. According to item
$13d_6)$ the last statement contradicts item $13a_5)$.

$13d_9)$. From item $13d_7)$ it follows that all results of
Chapter 2 are false.

$13e)$. Chapter 3 investigates a topological non-discrete free
algebra \break $F(X, K)$ of quasivariety $K$ in sense of item
$13c_1)$ with defining finite space $X$. In such a form all
results are false.

According to item 9g) let us transform $K$ into an almost
quasivariety. Then by Theorem 7 a free algebra $F(X, K)$ exists.
Since $X$ is a finite space then by item $9h_1)$ a non-discrete
free algebra $F(X, K)$ does not exist.

$13e_1)$. This confirmed again that all results of Chapter 3 are
false, are senseless.

$13f)$. We pass to Chapter 4. Literally we quote the items
$13f_1)$ -- $13f_6)$.

$13f_1)$. \textmd{Let $E = \cup\{E_n: n \in N_1 = \{0, 1\}$. The
set $E$ consists of the null-ary operations, and $E_1$ is the set
of one-ary operations $\ldots$ Fix a complete $T_i$-variety $K$ of
one-ary topological $E$-algebras.}

$13f_2)$. \textmd{Let $P$ be the set of all polynomials $\ldots$
It is clear that $P$ consists of all null-ary polynomials $P_0$
and one-ary polynomials $P_1$ (i.e. $P = P_0 \cup P_1$). If $p \in
P_0$, then for every $G \in K$, $p$ determines the element $e_p$.}

$13f_3)$. \textmd{\textbf{Definition 4.1.1.} Two polynomials $p$
and $q$ are equivalent, $p \sim q$, if:}

\textmd{1.$p(x) = g(x)$ is an identity, for $p, q \in P$. }
\textmd{2.  $p(x) = e_q$ is an identity, for $p \in P_1$, $q \in
P_0$.}

\textmd{3. $e_p  = e_q$  is an identity, for $p, q \in P _0$.}

$13f_4)$. \textmd{Obviously, the relation $p \sim q$ is an
equivalence relation on $P$. Fix $Q \subset P$ such that $Q$
contains only one element in every class of equivalence. Let $e
\in Q$, where $e(x) =x$ for all $x \in G \in K$. If $F(X, K)$ is
one-ary free algebra and $X \subset F(X, K)$, then $F(X,K) =
\cup\{p(X): p \in Q\} = \cup\{p(X): p \in P\}$.}

$13f_5)$. \textmd{Fix a discrete signature $E = E_0 \cup E_1$. Let
$K$ be a complete $T_i$-variety of topological $E$-algebras, where
$-1 \leq i \leq 3,5$.}

\textmd{\textbf{Theorem 4.2.1. a)} There are such discrete spaces
$D_0$ and $D_1$ so that for each $T_i$-space $X$ the topological
spaces of $F(X, K)$ and $D_0 \oplus (X \times D_1)$ are
homeomorphic.}

\textmd{\textbf{b} If $X$ is paracompact, or weakly paracoompact,
or metrizable, or $k$-space, then $F(X, K)$ is of the same type.}

\textmd{\textbf{c}  $\text{dim}$ $X$ = $\text{dim}$ $F(X, K)$,
$\text{ind}$ $X = $\text{ind}$ $F(X, K)$, $\text{Ind}$ $X =
$\text{Ind}$ $F(X, K)$.}

$13f_6)$.\textmd{\textit{Proof.} Let $Q \subseteq P(E)$ be as in
the Section 4.1 $\ldots$ The proof is complete.}

$13f_7)$. \textmd{\textbf{Theorem 4.2.2.} If $F(X, K)$ and $F(Y,
K)$ are free topologically one-ary algebras of the spaces $X$ and
$Y$ in the class $K$, then the following statements are
equivalent:}

\textmd{1. $X$ and $Y$ are homeomorphic.}

\textmd{$F(X, K)$ and $F(Y, K)$ are  topologically isomorphic.}

$13f_8)$. For item $13f_1)$ see the item $13c_2)$, i.e. a complete
$T_i$-variety of algebras not exist.

$13f_9)$. The definition of polynomial see the item $4c_4)$. By
item $13f_1)$ $K$  is a complete variety of one-ary algebras.
Hence $E = E_1$ and $E_0 = \{\emptyset\}$. By items $13f_2)$,
$13f_3)$ $P = P_0 \cup P_1$) and $P_0 \neq \{\emptyset\}$. If $p
\in P_0$ then from item $13f_2)$ it follows that  $p$ determines
the element $e_p$ in every $G \in K$. Hence $p \in E_0$.
Contradiction. Hence $P_0 = \{\emptyset\}$. In such case the
conditions 2, 3 from $13f_3)$ are without sense.

$13f_{10})$. The assertions from item $13f_4)$ are also absolutely
senseless. It seems that the author of monograph \cite{Chir1} does
not even know the definition of equivalence relation on set. For
more information regarding the absurdity of assertions from item
$13f_4)$ see the items $13a_1)$ -- $13a_6)$. The senseless of
items $13f_3)$ and $13f_4)$ also follows from the following item.

$13f_{11})$. According to items $13f_8)$, $13f_9)$ we correct
items $13f_1)$, $13f_3)$. We present \cite[pag. 348 -- 356]{Mal1}.
We consider the signature $\Omega = \{f_i \vert i \in I\}$, where
$f_i$ is a symbol of one-ary operation, i.e. $\Omega = \Omega_1$.
We denote by $V = \{v_i \vert i \in I\}$ a set of defining symbols
in class of semigroups. Let $K$ be the variety of
$\Omega$-algebras (in sense of Proposition 2) defined by
identities $f_{i_{\lambda_ 1}} \ldots
f_{i_{\lambda_{s\lambda}}}(x) = f_{j_{\lambda_ 1}} \ldots
f_{j_{\lambda_{t\lambda}}}(x) \quad (f_i, f_j \in \Omega)$, and
let $\overline{K}$ be the variety of semigroups (without unity)
defined by relations $v_{i_{\lambda_ 1}} \ldots
v_{i_{\lambda_{s\lambda}}} = v_{j_{\lambda_ 1}} \ldots
v_{j_{\lambda_{t\lambda}}}$.

By \cite[Theorem 1]{Mal1} an identity $f_{i_1} \ldots f_{i_s}(x) =
f_{j_1} \ldots f_{j_t}(x)$ is valid in variety $K$ if and only if
the corresponding identity $v_{i_1} \ldots v_{i_s} = v_{j_1}
\ldots v_{j_t}$ is valid in variety $\overline{K}$.

$13f_{12})$. The proofs of theorems from items $13f_5)$ and
$13f_7)$ are false. They use the false assertions from items
$13f_3)$ -- $13f_5)$ (see item $13f_6)$). Moreover, it ignores the
condition $X \subset F(X, K)$ from item $13f_4)$.

$13f_{13})$. According to the afore-mentioned we find that all
results from Chapter 4 are false, the Chapter 4 is a complete
senseless.

Now let us present the beginning of Chapter 5.

$13g)$. \textmd{In this section we examine the solvability of
equations over universal and free universal algebras. We also
present methods of constructing some free objects. The
construction of a free object in the class V ('E) and the results
from section 1.3.5 were also obtained by M.M. Choban using other
techniques \cite{ChobKir54}. The main result establishes that each
of the equations $ax = b$ and $ya = b$ over a free primitive
gruppoid has no more than two solutions.  The results of this
section were published in \cite{Kir133}.}

$13g_1)$. All results from Chapter 5 and work \cite{ChobKir54} are
totally false, are senseless.

$13g_2)$. Esteemed authors, please compare you researches with
items $13a_1)$ -- $13a_5)$, compare with Proposition 13, compare
with \cite[Chapter1]{Smith1}, where the construction of free
groupoids with division, of free quasigroups is described in
detail.

$13h)$. We pass to the brief analysis of Chapter 7, whose results
are published in  \cite{Chir135}. We ascertain that all results
from Chapter 7  are totally false, are senseless. We confirm it.

$13h_1)$. \textmd{\textbf{Theorem 7.2.2} The free groupoid with
division $\Gamma(X)$ is a quasigroup.}

$13h_2)$. In $13h_1)$ $\Gamma(X)$ is defined as  a groupoid for
which the equations $ax = b$ and $ya = b$ ($\forall a, b \in
\Gamma (X)$) have solutions. In such a case $\Gamma(X)$ is not an
universal algebra.

$13h_3)$. From item $13h_2)$ it follows that the Theorem 7.2.2 is
false. From here it follows that the Definitions 5.5.3 and 7.2.1
of M. Choban, the Theorem 5.5.5 and 7.2.4, the Corollary 7.2.3 are
false.

$13i)$. Let us analyze Chapter 8. We quote items $13i_1)$,
$13i_4)$,  $13i_7)$, $13i_8)$.

$13i_1)$. \textmd{A $T_2$-space $X$ is called a $k$-space if the
subset $F \subset X$ is closed in $X$ iff $F \cap \Phi$ is compact
for each compact $\Phi \subset X$} (pag. 22).

$13i_2)$. The definition of $T_i$-quasivariety of $k$-$E$-algebras
follows from item 9c)  if instead of topological $E$-algebras we
consider a topological $k$-$E$-algebras (for comments see item
9g).)

$13i_3)$. The definitions of complete $T_i$-quasivariety from
Chapter 1 (item $13c_1)$) and Chapter 8 are totally different.

$13i_4)$. \textmd{ A mapping $f: X \rightarrow Y$ of a space $X$
into a space $Y$ is called a $k$-continuous mapping if for every
compact subset $\Phi \subseteq X$ the restriction $f|\Phi: \Phi
\rightarrow Y$ is continuous} (pag. 103).

$13i_5)$. For a $T_i$-quasivariety of $k$-$E$-algebras $V$ The
definitions of free algebra $(F^a(X, V), a_X)$, of $k$free algebra
$(F^k(X, V), q_X)$, of $t$-free algebra $(F(X, V), t_X)$ of the
space $X$ follow from items 10g) and 10h) or 6c), if instead of
continuous mappings we consider $k$-continuous mappings (for
comments see items $10g_1)$, $10h_1)$, $10h_2)$.)

$13i_6)$. For the definitions from items $13i_2)$, $13i_5)$ it is
necessary that the signature $E$ of $E$-algebras be a $T_i$-space,
i.e. for the $E$ to be a topological space. According to item 1e)
this requirement is a fiction.

$13i_7)$. \textmd{\textbf{Theorem 8.3.4.} Let $V$ be a
$T_i$-quasivariety of $k-E$-algebras. Then for every non-empty
space $X$ there exist:}

\textmd{1.  the unique free $E$-algebra $(F^a(X; V); a_X)$;}

\textmd{2, the unique $k$-free $E$-algebra $(F^k(X; V); q_X)$;}

\textmd{3. the unique $t$-free $E$-algebra $(F(X; V ; t_X)$;}

 \textmd{4.  the unique continuous homomorphisms
$b_X: F^a(X; V ) \rightarrow F^k(X; V)$, $c_X : F^a(X; V )
\rightarrow F(X; V)$ and $l_X: F^k(X; V) \rightarrow F(X; V)$ such
that $q_X = b_X\circ a_X$ and $t_X = c_X \circ a_X = l_X \circ
q_X$.}

$13i_8)$. \textmd{\textit{Proof.} Let $\tau$  be an infinite
cardinal and $\tau \geq |X| + |E|$. Consider that $\{f_{\alpha}: X
\rightarrow G_{\alpha} \in  V: \alpha \in A\}$ is the set of all
mappings for which $|G_{\alpha}| \leq \tau$.
 Let $C = \{\alpha \in A:
f_{\alpha}: f_{\alpha}: X \rightarrow G_{\alpha} \text{is
continuous} \}$ and $D = \{\alpha \in A: f_{\alpha}: f_{\alpha}:
 X \rightarrow G_{\alpha} \text{is k-continuous}
\}$. Fix a subset $B \subseteq A$. Consider the diagonal product
$f_B: X \rightarrow G_B = \prod\{G_{\alpha}: \alpha \in B\}$,
where $f_B(x) = \{f+{\alpha}(x): \alpha \in B\}$. Let $F(X; B)$ be
the subalgebra of $G_B$ generated by the set $g_B(X)$. Consider
the projection $g_{\alpha}: F(X, B) \rightarrow G_{\alpha}$ for
every $\alpha \in  B$.}

\textmd{If $\alpha \in B$, then $f_{\alpha} = g_{\alpha} \circ
f_B$. Fix a mapping $f: X \rightarrow G \in V$.}

\textmd{Put $G^{\prime} = a(E; f(X))$. Then $|G^{\prime}| \leq
\tau$ and $f = f_{\alpha}$, $G^{\prime}  = G_{\alpha}$ for some
$\alpha \in A$.}

\textmd{If $\alpha \in B$, then $f = g_{\alpha} \circ f_B$ and
$g_{\alpha}: F(X;B) \rightarrow G$ is a continuous homomorphism.
Hence $(F^a(X; V ); a_X) = (F(X; A); f_A)$, $(F^k(X; V ); q_X) =
(F(X;D); f_D)$ and $(F(X; V ); t_X) = (F(X;C); f_C)$. The
assertions are \break  proved.}

$13i_9)$. The proof from item $13i_8)$ almost literally repeats
the proof from item 6d, i.e. from \cite{Chob19}. See also the item
2c). The proof from item $13i_8)$ is false, is completely
senseless like the proof from item 6d. Foe comments see item
$6d_1)$.

All results from Chapter 8 use essentially Theorem 8.3.4 from item
$13i_7)$. Hence all results from Chapter 8 are not proved, more
specific, are senseless. This is also proved by item $13i_{10})$.

$13i_{10})$. In sections 8.4 -- 8.6, 8.8, 8.10 it is necessary
that signature $E$ of $E$-algebras be a $k_{\infty}$-signature,
i.e. that $E$ be a topological space. According to item 1e) this
requirement is a fiction (see also item $13i_6)$).

$13j)$. Let us analyze Chapter 6, where an uniform structures is
applied to study the free topological algebras. For the theory of
uniform spaces see, for example, \cite[Chapter 8]{Eng}. Chapter 6
investigates the notions from items $13j_1)$ -- $13j_5)$,
$13j_{11})$ $13j_{12})$. We state them.

$13j_1)$. \textmd{\textbf{Definition 6.5.1} A class $K$ of uniform
$E$-algebras presents a full variety, provided that:}

\textmd{1. the class $K$ is closed under Cartesian products and
under taking subalgebras;}

\textmd{2. the class $K$ is closed under passing to uniformly
continuous homomorphic images;}

\textmd{3. if $\mu A \in K$ and the algebra $A$ is uniform with
respect to a uniformity $\eta$, then $\eta A \in K$.}

$13j_2)$. \textmd{\textbf{Definition 6.6.3} A class $K$ of
$T$-uniformizable algebras is called a full $T$-uniformizable
variety, provided that:}

\textmd{1. the class $K$ s closed under Cartesian products and
taking subalgebras;}

\textmd{2. if $A \in K$ and $A$ is a $T$-uniform algebra with
respect to a topology $\tau$ and a uniformity $\mu$, then $(A;
\tau) \in K$;}

\textmd{3. if $A \in K$ and $f: A \rightarrow B$ is a continuous
homomorphism onto a $T$-uniformizable algebra $B$ then $B \in K$.}

$13j_3)$. \textmd{Let $K$ a nontrivial full $T$-uniformizable
variety of topological $E$-algebras. Then, for every algebra $A
\in K$, we can define the maximal uniformity $\nu_A$ that
$T$-uniformizes $A$.}

$13j_4)$. \textmd{\textbf{Theorem 6.6.4} Let $K$ be a nontrivial
variety of $E$-algebras. Then, for every uniform space $\mu X$,
there exists an algebra $F(\mu X; K) \in K$ such that:}

\textmd{1. $\mu X$ is a uniform subspace of the uniform space
$(F(\mu X; K); \nu_{F(\mu X;K)})$;}

\textmd{2. the set $X$ generates $F(\mu X; K)$ algebraically;}

\textmd{3. for every uniformly continuous mapping $f : X
\rightarrow A \in K$, where $A$ is endowed with the uniformity
$\nu_A$, there exists a continuous homomorphism $\overline{f} :
F(\mu X; K) \rightarrow A$ such that $f = \overline{f}|X$;}

\textmd{4. the algebra $F(\mu X; K)$ is algebraically free in the
class $K$.}

$13j_5)$. \textmd{\textit{Proof.} We denote by $K_1$ all the
algebras of $K$ endowed with uniformities that make them uniform.
Then $K_1$ is a full variety of uniform algebras. Each algebra $A
\in K$ with the discrete uniformity belongs to $K_1$. Therefore,
$K$ and $K_1$ coincide algebraically $\ldots$}

$13j_6)$. The definitions from items $13j_1)$ and $13j_2)$ are
false. According to item 9i) a nontrivial full variety of uniform
$E$-algebras and a nontrivial full $T$-uniformizable variety of
$T$-uniformizable algebras do not exist.

$13j_7)$. We mentioned that all results from Chapter 6 are proved
for nontrivial full varieties of topological $E$-algebras, and
their incorrectness follows from item $13j_6)$. Particularly, see
item $13j_5)$.

$13j_7)$. Moreover, item $13j_8)$ proves once more the
incorrectness and senselessness of papers \cite{AC14}, \cite{DC},
\cite{Chob3}, \cite{Chob191}.

$13j_8)$. \textmd{Let $K$ be a full nontrivial $T_i$-quasivariety
of topological $E$-algebras. The articles \cite{AC14}, \cite{DC}
prove that a free topological algebra always exists and is unique
to within a topological isomorphism. If $X$ is completely regular
then $i_X: X \rightarrow F(X;K)$ is a topological embedding
\cite{AC14}, \cite{DC}. If $X$ is discrete then $F(X;K) =
F^a(X;K)$ and $i_X = j_X$. Therefore, for any space $X$, the
algebra $F^a(X;K)$ exists, is unique, and presents a discrete
space. There always exists a continuous onto homomorphism $k_X :
F^a(X;K) \rightarrow F(X;K)$ such that $i_X = j_X \circ k_X$. If
$k_X$ is an isomorphism then the algebra $F(X;K)$ is referred to
as abstractly free. For any Tychonoff space $X$, the algebra
$F(X;K)$ is abstractly free \cite{Chob3}. The following fact plays
an important role in studying free objects \cite{Chob191}.}

$13j_{9})$. The content of Theorem 6.6.4 is floppy. In condition 1
it is used the maximal uniformity $\nu_A$, as defined for full
$T$-uniformizable variety of topological $E$-algebras according to
item $13j_3)$.

$13j_{10})$.  In item $13j_5)$ we presented only the false begin
of proof. Further the proof uses essentially condition 4 of the
Theorem. Incorrectness of condition 4 follows from item $13j_8)$.
Consequently, the Theorem 6.6.4 from Section 6.6 is false.

$13j_{11})$. All the following Sections 6.7 - 6.12 contain
something unbelievable, unexplainable, a real pun of topological
notions. We confirm it.

$13j_{11})$. \textmd{A topological algebra $A$ is $T$-uniformizale
if on $A$ there exists a compatible uniformity that induces the
structure of a $T$-uniform algebra on $A$} (pag. 80).

$13j_{12})$. \textmd{Let us fix a signature $E$, $i \in \{-1; 0;
1; 2; 3; 3,5\}$, and a nontrivial full $T_i$-variety $K$ of}
topological $E$-algebras. We denote by $K_u$ the collection of all
$T$-uniformizable algebras of $K$. Clearly, $K_u$ contains all
paracompact, all Diedonne-complete, and all discrete algebras of
$K$. \textmd{Thus, $K_u$ is a full T-uniform variety of algebras
and the classes $K$ and $K_u$ coincide algebraically.
 Hence, for every Tychonoff space $X$, there exists a
continuous isomorphism $\pi_X:  F(X; Y) \rightarrow F(X; K_u)$
such that $\pi_X = x$ for all $x \in X$ and $\pi_X|X$ is a
homeomorphism} (pag. 82)

$13j_{13})$. One of the main results of Section 6.7 - 6.12
(Theorems 6.9.2, 6.9.3, 6.10.6, 6.11.2, Corollary 6.12.1) are the
statements of form:

Let $i = 3,5$. For a Tychonoff space $X$, the following conditions
are equivalent:

$X$ is a $\prod_{\omega}$-space (respect. $C_{\omega}$-space, or
Diedonne-complete space, or Lindelof $P$-space),

2. $F(X; K)$ is a $\prod_{\omega}$-space (respect.
$C_{\omega}$-space, or Diedonne-complete space, or Lindelof
$P$-space),

3. $F(X; K_u)$ is a $\prod_{\omega}$-space (respect.
$C_{\omega}$-space, or Diedonne-complete \break space, or Lindelof
$P$-space).

Similar statement holds in case when $K$ is a nontrivial full
variety of topological groups with operators.

$13j_{14})$. The incorrectness of assertions from items
$13j_{12})$. $13j_{13})$ is indicated in item $13j_6)$. Moreover,
the assertions from items $13j_{12})$. $13j_{13})$ contradict
Theorem 7 from the first part of present paper, $K$ should be an
almost quasivariety of topological algebras. In this case the
space $X$ should meet conditions $K_1)$, $K_2)$ from Section 2,
but by no means the topological space of free algebra $F(X; K)$.
Hence the equivalence of conditions 1, 2 from item $13j_{13})$ is
senseless. It seems that the author of monograph is unaware of the
most basic notions: does not understand the sense and meaning of
the notion of topological algebra given by defining space
(Definition 1 of Section 2 of first part or the definition of
topological free algebra of a space from item $13c_3)$. But
compare with the proofs offered by the author for the assertions
from item $13j_{13})$.

$13j_{15})$. According to items $13j_{6})$, $13j_{7})$,
$13j_{9})$, $13j_{10})$, $13j_{11})$,  $13j_{14})$ we conclude
that Chapter 6, like Chapter 8 (see item $13i_9)$), is
anti-scientific. Publishing such chapters is a crime.

$13j_{16})$. On pages 68 and 102 it is mentioned that Chapters 6
and 8 were published in \cite{ChobKir4} and \cite{ChobKir5}.  We
mentioned that Chapters 6 and 8 coincide with Chapters 2 and 4 of
Thesis for a Habilitat Doctors Degree \cite{Chir2} and form its
foundation.

$13k)$. Let us analyze now Chapter 7. It is stated literally.

$13k_1)$. \textmd{An $n$-groupoid $A$ with $n$-ary operation
$\omega$ is called an $n$-groupoid with division or an
$nD$-groupoid, if the equation $\omega(a_1, \ldots , a_{i-1}, x,
a_{i+1}, \ldots \break \ldots , a_n)  = b$ has a solution (not
necessarily unique), for every $a_1, \ldots , a_n, b \in A$  and
any $1 \leq i \leq n$.}

$13k_2)$. If in item $13k_1)$ $n = 2$, then $(A, \cdot)$ is called
a groupoid with division. This means that the equations $a\cdot x
= b, y\cdot a = b$ have a solution, for every $a, b \in A$.

$13k_3)$. \textmd{Define a free object of a set $X$ in the class
$V(E)$ of all groupoids with division according to next definition
of M. Choban.}

\textmd{\textbf{Definition 3.13.1.} The free groupoid with
division of a set $X$ in the class $V(E)$,  is an E-algebra
$\Gamma(X) \in V(E)$ such that:}

\textmd{1. $X \subset \Gamma(X)$ and the set $X$ algebraically
generates the algebra $\Gamma(X)$, i.e. if $X \subset Y \subset
\Gamma(X)$, $Y \neq \Gamma(X)$, and $Y$ is a subalgebra of the
algebra $\Gamma(X)$, then $Y \notin V(E)$.}

\textmd{2. For every mapping $f  X \rightarrow A$, where $A \in
V(E)$, there exists a homomorphism $\overline{f}: \Gamma(X)
\rightarrow A$ such that $\overline{f}|X = f$.} 117

$13k_4)$. \textmd{\textbf{Theorem 3.13.2.} The free groupoid with
divisions $\Gamma(X)$ is a quasigroup.}

$13k_5)$. \textmd{\textit{Proof.} It follows from Theorem 5.5.5.
Nevertheless, we will present a direct proof.}

$13k_6)$. The $n$-groupoid with division from item $13k_1)$ is not
an universal algebra.

$13k_7)$. According to item $13k_6)$ the definition from item
$13k_3)$ is false. It contradicts Theorem 7 from the first part,
$V(E)$ should be an almost quasivariety of universal algebras.

$13k_8)$. According to item $13k_6)$ the theorem  from item
$13k_3)$ is false. Moreover, if defining the groupoid with
divisions algebraical, i.e. using three binary operations, then
the theorem is false anyway.

$13k_9)$. From item $13k_5)$ it follows that the Theorem 5.5.5 is
false. Clearly, the presented proof is false, is cumbersome and
unclear. .

$13k_{10})$. Theorem 7.2.4 gives a conditions when continuous
homomorphisms of topological groupoids with continuous
homomorphisms of topological \break groupoids with a continuous
division are open.

$13k_{11})$. The theorem from item $13k_{10})$ is false. Its proof
uses the false theorem from item $13k_8)$ and the Mal'cex's
Theorem \cite{Mal2} that the topological quasigroup is regular. A
necessary condition that a continuous homomorphism of groupoid
with division is given in Corollary 21 and sufficient condition is
given in Corollary 2 from the first part or \cite[Theorem
10]{Mal2}.

$13k_{12})$. From items $13k_{6})$ -- $13k_{11})$ it follows that
all results from Chapter 7 are false.

$13l)$. On page 125 it is mentioned that the results from Chapter
9 were published in \cite{ChobChir55}. We present the items
$13i_1$, -- $13i_3$ from the beginning of this chapter.

$13l_1)$. \textmd{L. S. Pontrjagin \cite{Pont} proved that a
linear connected space that covers a topological group admits, in
a natural way, a structure of a topological group. In this work we
establish a similar result for universal algebras with continuous
signature. This result, for the case of a finite discrete
signature, was obtained by A.I. Mal'cev \cite{Mal2}.}

$13l_2)$. \textmd{Result from this work is stronger then Mal'cev's
Theorem. In particular, the result holds for the topological
R-modules, where R is a topological ring.}

$13l_3)$. \textmd{We mention that if $J$ is a set of identities,
then the totality $V(J)$ of all Hausdorff topological
$E$-algebras, which satisfy the identities $J$, forms a complete
variety of topological $E$-algebras.}

$13l_4)$. According to item 1e) the property of continuous
signature should be ignored. Hence the result from this work not
is stronger then Mal'cev's Theorem. In particular, the assertion
that the result holds for the topological R-modules, where R is a
topological ring, is false.

$13l_5)$. The assertion from item $13i_3)$ is false.

$13l_6)$. In Chapter 9 a complete varieties of algebras are
investigated. According to item $13c_2)$ a complete varieties does
not exist.

$13l_7)$. The items $13l_4)$ -- $13l_6)$ are enough to state:
Chapter 9 is senseless.

$13m)$. In Chapters 10 -- 14 some notions are investigated, that
present little algebraic interest: medial quasigroups with left
unity, paramedial quasigroups. multiple identities, Fuzzy algebras
and others. The veracity of their results is not different than
the veracity of the results from the previous Chapters. For an
example, we will analyze briefly Chapter 14.

$13n)$. This section provides a study of the category of fuzzy
groupoids with division. The results of this section were
published in \cite{ChobChir58}, \cite{ChobChir60}.

$13n_1)$. \textmd{\textbf{Definition 6.7.3} An $E$-algebra $G$ is
an:}

\textmd{1. $E$-groupoid with left division if there exist the
operations $A; G; W \in E_2$ for which $A(C(y; x); x) = W(C(x; y);
x) = y$ for every points $x; y \in G$.}

\textmd{2. $E$-groupoid with right division if there exist the
operations $A; B; V \in E_2$ for which $A(x;B(x; y)) = V (x;B(y;
x)) = y$ for every $x; y \in G$.}

\textmd{3. $E$-groupoid with division if there exist the operation
$A; B; G; V; W \in E_2$ such that $A(x; B(x; y)) = A(G(y; x); x) =
V (x; B(y; x)) = W(G(x; y); x) \break = y$ for every $x; y \in
G$.}

\textmd{4. $E$-quasigroup if there exist the operations $A; B; C
\in  E_2$ for which $A(x; B(x; y)) = B(x; A(x; y)) = A(G(y; x); x)
= G(A(y; x); x) = y$.}

\textmd{Every $E$-groupoid with division is an $E$-groupoid with
left and right divisions. If $G$ is an $E$-groupoid, then $G$ is
also an $E$-groupoid with division and $V = G; W = B$.}

$13n_2)$. I believe it is not even worth commenting it, taking
into account at least the identities (20) -- (22) from the first
part.

From the above-mentioned we conclude that the results from
monograph \cite{Chir1} and the thesis for a Habilitat Doctors
Degree \cite{Chir2} are false, do not meet the minimum
requirements for such works.

Esteemed authors of paper \cite{ChobChir44}, I think that it is
unbelievable that such works were published. I believe it is a
unique phenomenon in this world. To fool the readers and the
journals, where these works were published, for more than 40 years
is a real crime.

Nicolae I. Sandu,

Tiraspol State University of Moldova,

Chisin\u{a}u, R. Moldova

sandumn@yahoo.com

\end{document}